\documentclass[11pt]{article}

\usepackage{amsfonts}
\usepackage{amssymb}
\usepackage{amsmath}
\usepackage{amsthm}

\newcommand{\nicefrac}[2]
{\leavevmode \kern.1em\raise.5ex\hbox{\the\scriptfont0 #1}
             \kern-.1em/\kern-.15em\lower.25ex
             \hbox{\the\scriptfont0 #2}}

\newtheorem*{theorem}{Theorem}
\newtheorem*{proposition}{Proposition}
\newtheorem*{corollary}{Corollary}
\newtheorem*{problem}{Problem}
\newtheorem*{definition}{Definition}

\newtheorem*{remark}{Remark}
\newtheorem*{remarks}{Remarks}

\theoremstyle{definition}

\begin{document} 

\begin{center}
{\Large{\sc On $2$-surfaces in $\mathbb R^4$ and $\mathbb R^n$}}\\[1cm]
{\large Steffen Fr\"ohlich}\\[0.4cm]
{\small\bf Abstract}\\[0.4cm]
\begin{minipage}[c][2.5cm][l]{12cm}
{\small In this overview report we generalize Erhard Heinz' curvature estimate for minimal graphs in $\mathbb R^3$ to graphs in $\mathbb R^n$ of prescribed mean curvature. Secondly, we analyse these problems in the frame of the outer differential geometry which leads us to the notions of normal torsion and normal curvature for immersions in $\mathbb R^4.$}
\end{minipage}
\end{center}
{\small MCS 2000: 35J60, 53A07, 53A10}
\subsection{Immersions in $\mathbb R^n.$ Basic notations}
On the closed unit disc $B:=\{(u,v)\in\mathbb R^2\,:\,u^2+v^2\le 1\}$ we consider two-dimensional immersions
\begin{equation}\label{1.1}
  X(u,v)=\big(x^1(u,v),x^2(u,v),\ldots,x^n(u,v)\big),\quad(u,v)\in B,
\end{equation}
of class $X\in C^3(B,\mathbb R^n),$ $n\ge 3.$\\[1ex]
Let $u^1\equiv u$ and $u^2\equiv v,$ and denote by $X_{u^i},$ $X_{u^i}^t$ the partial derivative of $X$ and its transpose w.r.t. $u^i.$ We define the coefficients of the first fundamental form of $X$ as
\begin{equation}\label{1.2}
  g_{ij}:=X_{u^i}\cdot X_{u^j}^t\,,\quad i,j=1,2.
\end{equation}
The surface area element $W$ of an immersion fulfills
\begin{equation}\label{1.3}
  W:=\sqrt{g_{11}g_{22}-g_{12}^2}>0\quad\mbox{in}\ B.
\end{equation}
We preferably use conformal parameters $(u,v)\in B$ with the properties
\begin{equation}\label{1.4}
  |X_u|^2=W=|X_v|^2\,,\quad
  X_u\cdot X_v^t=0
  \quad\mbox{in}\ B.
\end{equation}
With such an immersion we associate a regular frame
\begin{equation}\label{1.5}
  \Big\{X_u,X_v,N_1,\ldots,N_{n-2}\Big\}
\end{equation}
which consists of the two tangential vectors $X_u,$ $X_v,$ and an orthonormal (ON) system $N_\Sigma,$ $\Sigma=1,\ldots,n-2,$ in the normal space of $X.$
\subsection{Minimal surface graphs in $\mathbb R^3$}
\setcounter{equation}{0}
\begin{theorem}
(S.\,Bernstein, 1916; E.\,Heinz, 1952)\\[0.1cm]
Let the minimal surface graph $X(x,y)=(x,y,\varphi(x,y)),$ where $\varphi\in C^3(B_R,\mathbb R)$ on $B_R:=\{(x,y)\in\mathbb R^2\,:\,x^2+y^2\le R^2\}$ solves the non-parametric minimal surface equation
\begin{equation}\label{2.1}
  (1+\varphi_y^2)\varphi_{xx}-2\varphi_x\varphi_y\varphi_{xy}+(1+\varphi_x^2)\varphi_{yy}=0
  \quad\mbox{in}\ B_R\,,
\end{equation}
be parametrized using conformal parameters $(u,v)\in B$ as a vector mapping $X=X(u,v).$ Then there is a universal constant $\Theta\in[0,+\infty)$ such that
\begin{equation}\label{2.2}
  \kappa_1(0,0)^2+\kappa_2(0,0)^2\le\frac{\Theta}{R^2}
\end{equation}
holds true for its principle curvatures $\kappa_1$ and $\kappa_2$ at the origin $(0,0)\in B.$
\end{theorem}
\renewcommand{\proofname}{Idea of the proof}
\begin{proof}
The unit normal vector of the immersion is given by
\begin{equation}\label{2.3}
  N(u,v):=\frac{X_u(u,v)\times X_v(u,v)}{|X_u(u,v)\times X_v(u,v)|}\,.
\end{equation}
Using conformal parameters $(u,v)\in B,$ it holds the estimate
\begin{equation}\label{2.4}
  \kappa_1(0,0)^2+\kappa_2(0,0)^2
  \le\frac{|\nabla N(0,0)|^2}{W(0,0)}\,.
\end{equation}
From $|\triangle N|\le|\nabla N|^2,$ where $\triangle$ and $\nabla$ mean the usual Laplacian and gradient, together with a suitable modulus of continuity for $N,$ one deduces an universal constant $C_1\in(0,+\infty)$ such that $|\nabla N(0,0)|\le C_1.$ Secondly, $\triangle X=0$ implies the lower estimate $W(0,0)\ge C_2R^2,$ $C_2\in(0,+\infty).$ Now set $\Theta:=\frac{C_1^2}{C_2}.$
\end{proof}
\begin{remark}
The  proof is based on a detailed analysis of non-linear elliptic systems with quadratic growth in the gradient
\begin{equation}\label{2.5}
  |\triangle X|\le a|\nabla X|^2\quad\mbox{in}\ B,\quad\mbox{where}\ a\in\mathbb R_+\,.
\end{equation}
\begin{itemize}
\item[(i)]
If $a\cdot\sup|X|<1,$ then there is a constant $c_1=c_1(a,M)$ such that
\begin{equation}\label{2.6}
  |\nabla X(0,0)|\le c_1(a,M),\quad
  M:=\sup|X(u,v)|.
\end{equation}
\item[(ii)]
For plane mappings $z\colon\overline B\to\mathbb R^2,$ with $z\colon\partial B\to\partial B$ topologically and positively oriented, $z(0,0)=(0,0),$ $J_z(w)>0$ in $B$ for its Jacobian, and which satisfy (\ref{2.5}) with $0<a<\frac{1}{2},$ there is a constant $c_2=c_2(a)$ with
\begin{equation}\label{2.7}
  |\nabla z(0,0)|\ge c_2(a).
\end{equation}
\end{itemize}
Analog estimates can be proved on small inner discs $B_\varrho(0,0),$ $\varrho\le\varrho_0<1.$ We refer to \cite{Heinz_01} and \cite{Sauvigny_01}.
\end{remark}
\noindent
Letting $R\to\infty$ in (\ref{2.2}) yields the
\begin{corollary}
A complete minimal graph of regularity class $C^3$ is a plane.\footnote{In the proof we applied the Laplacian $\triangle$ to $N$ which assumes $C^3$-regularity. In fact, $C^2$ is sufficient.}
\end{corollary}
\subsection{Minimal surface graphs in $\mathbb R^n$ for $n\ge 3$}
\setcounter{equation}{0}
We apply the above method to graphs in higher codimensional Euclidean spaces.
\begin{theorem}
Let the minimal graph
\begin{equation}\label{3.1}
  X(x,y)=(x,y,\varphi_1(x,y),\ldots,\varphi_{n-2}(x,y)),\quad\varphi_\Sigma\in C^3(B_R,\mathbb R),
\end{equation}
be parametrized in conformal parameters $(u,v)\in B.$ Then there is a universal constant $\Theta\in[0,+\infty)$ such that for its Gaussian curvature $K_N$ along any unit normal section $N$ of the frame (\ref{1.5}) it holds
\begin{equation}\label{3.2}
  |K_N(0,0)|\le\frac{\Theta}{R^4}\,\|X\|_{C^0(B_R)}^2\,,\quad
  \|X\|_{C^0(B_R)}:=\sup_{(x,y)\in B_R}|X(x,y)|.
\end{equation}
\end{theorem}
\begin{proof}
Using conformal parameters $(u,v)\in B,$ we estimate
\begin{equation}\label{3.3}
  |K_N(0,0)|\le\frac{|X_{uu}(0,0)||X_{vv}(0,0)|+|X_{uv}(0,0)|^2}{W(0,0)^2}\,.
\end{equation}
Due to $\triangle X=0$ in $B,$ potential theory ensures a third constant $C_3\in(0,+\infty)$ with the property
\begin{equation}\label{3.4}
  |X_{u^iu^j}(0,0)|\le C_3\|X\|_{C^0(B_R)}
\end{equation}
for $i,j=1,2$ (see e.g. \cite{Gilbarg_Trudinger_01}, Theorem 4.6). As in the previous paragraph we deduce a lower bound for $W(0,0).$ Now, set $\Theta:=\frac{2C_3^2}{C_2^2}.$
\end{proof}
\begin{corollary}
(Bernstein-Liouville theorem)\\[0.1cm]
Assume that the minimal graph $X(x,y)$ satisfies the growth condition
\begin{equation}\label{3.5}
  \|X\|_{C^0(B_R)}\le\Omega R^\varepsilon,\quad
  \varepsilon\in(0,2),
\end{equation}
with a constant $\Omega\in(0,+\infty).$ Then, if $X$ is complete, it represents a plane.
\end{corollary}
\noindent
Namely, note that with (\ref{3.5}) it holds
\begin{equation}\label{3.6}
  |K_N(0,0)|\le\frac{\Theta\Omega^2}{R^4}\,R^{2\varepsilon}\longrightarrow 0
  \quad\mbox{for}\ R\to\infty
\end{equation}
for all $N,$ and the Corollary follows.
\begin{remarks}
\begin{itemize}
\item[1.]
The Corollary is sharp in the following sense: $(w,w^2)$ for $w\in\mathbb C$ is a non-plane minimal surface graph over $\mathbb C.$
\item[2.]
It generalizes the Liouville theorem in complex analysis because a
holomorphic function $\varphi(z)$ generates a minimal graph $(z,\varphi(z))$
in $\mathbb R^4.$
\end{itemize}
\end{remarks}
\subsection{Surfaces of prescribed mean curvature in $\mathbb R^4$}
\setcounter{equation}{0}
Given a vector field ${\mathcal H}\,:\,\mathbb R^4\to\mathbb R^4,$ define the scalar $H(X,Z):={\mathcal H}(X)\cdot Z^t$ for $X\in\mathbb R^4,$ $Z\in S^3:=\{X\in\mathbb R^4\,:\,|X|^2=1\}.$
\begin{definition}
The immersion $X\,:\,B\to\mathbb R^4$ is called a {\it conformally parametrized surface of prescribed mean curvature field} $H\,:\,\mathbb R^4\times S^3\to\mathbb R$ iff
\begin{equation}\label{4.1}
\begin{array}{l}
  \triangle X=2H(X,N_1)WN_1+2H(X,N_2)WN_2\,, \\[0.2cm]
  |X_u|^2=W=|X_v|^2\,,\quad X_u\cdot X_v^t\quad\mbox{in}\ B
\end{array}
\end{equation}
for an arbitrary ON-normal section $\{N_1,N_2\}.$
\end{definition}
\begin{remarks}
\begin{itemize}
\item[1.]
Note that the mean curvature field depends naturally on the space variable $X$ and the attached direction $Z.$ In this sense, the mean curvature field is Finslerian.
\item[2.]
If $X$ together with an ON-normal section $\{N_1,N_2\}$ solves (\ref{4.1}), then
\begin{equation}\label{4.2}
  |\triangle X|\le 2h_0|\nabla X|^2\quad
  \mbox{with}\ h_0:=\sup_{(X,Z)\in\mathbb R^4\times S^3}|H(X,Z)|.
\end{equation}
Thus, the immersion satisfies the structur condition (\ref{2.5}).
\end{itemize}
\end{remarks}
\noindent
Instead of applying potential theory to immersions (\ref{4.1}) directly (but see (\ref{4.7}) below), we refer to \cite{Bergner_Froehlich_01} for the following curvature estimate:
\begin{theorem}
Let $X\in C^{2+\alpha}(B,\mathbb R^4)$ be a surface of prescribed mean curvature field $H=H(X,Z)$ such that
\begin{itemize}
\item[(A1)]
$X$ is a positively oriented conformal reparametrization of a graph
\begin{equation}\label{4.3}
  X(x,y)=(x,y,\varphi_1(x,y),\varphi_2(x,y)),\quad\varphi_\Sigma\in C^{2+\alpha}(B_R,\mathbb R);
\end{equation}
\item[(A2)]
there hold $|H(X,Z)|\le h_0$ for all $X\in\mathbb R^4,$ $Z\in S^3,$ and
\begin{equation}\label{4.4}
  |H(X_1,Z_1)-H(X_2,Z_2)|\leq h_1|X_1-X_2|^\alpha+h_2|Z_1-Z_2|
\end{equation}
with real constants $h_0,h_1,h_2\in [0,+\infty);$
\item[(A3)]
$X$ represents a geodesic disc of radius $r>0$ such that
\begin{equation}\label{4.5}
  \int\hspace{-0.25cm}\int\limits_{\hspace{-0.3cm}B}|\nabla X|^2\,dudv\le d_0r^2
  \quad\mbox{with a constant}\ d_0>0;
\end{equation}
\item[(A4)]
At every point $w\in B,$ each normal vector of the immersion makes an angle of at least $\omega>0$ with the $x_1$-axis.
\end{itemize}
Then there is a constant $\Theta(h_0r,h_1r^{1+\alpha},h_2r,d_0,\sin\omega,\alpha)\in(0,+\infty)$ such that
\begin{equation}\label{4.6}
  \kappa_{\Sigma,1}(0,0)^2+\kappa_{\Sigma,2}(0,0)^2
    \le\frac{1}{r^2}\Big\{(h_0r)^2+\Theta\Big\}.
\end{equation}
\end{theorem}
\begin{proof}
One evaluates (\ref{3.3}). The graph property together with (A3) and (A4) ensure a variant of (\ref{2.7}) which estimates $W(0,0)$ from below. Next, bounding $W$ from above gives the Poisson problems
\begin{equation}\label{4.7}
  \triangle x^i=2H_1WN_1^i+2H_2WN_2^i=:f^i\,,\quad
  H_\Sigma:=H(X,N_\Sigma),\quad i=1,\ldots,4,
\end{equation}
with $\|f\|_{C^\alpha(B_\varrho)}\le C_5$ with a suitable constant $C_5\in(0,+\infty),$ $\varrho\in(0,1).$ Potential theory yields again upper bounds for $|X_{u^iu^j}(0,0)|.$
\end{proof}
\begin{remarks}
\begin{itemize}
\item[1.]
Assumption (A4) is a special case of a so-called Osserman-condition from \cite{Osserman_01}. The $x_1$-axis is chosen arbitrarily.
\item[2.]
For $H\equiv 0,$ the limit $r\to\infty$ gives a Bernstein type theorem for minimal surface graphs.
\end{itemize}
\end{remarks}
\subsection{Outer differential geometry in $\mathbb R^4$}
\setcounter{equation}{0}
\subsubsection{The Ricci equations and the normal curvature}
Let the immersion $X\in C^3(B,\mathbb R^4)$ be equipped with a regular moving frame $\{X_u,X_v,N_1,N_2\},$ where $\{N_1,N_2\}$ forms an ON-basis of the normal space.\\[1ex]
Denote by
\begin{equation}\label{5.1}
  L_{\Sigma,ij}:=X_{u^iu^j}\cdot N_\Sigma^t\,,\quad i,j=1,2,\ \Sigma=1,2,
\end{equation}
the coefficients of the second fundamental form w.r.t. $N_\Sigma,$ and by
\begin{equation}\label{5.2}
  \sigma_{\Sigma,i}^\Omega:=N_{\Sigma,u^i}\cdot N_\Omega^t\,,\quad i=1,2,\ \Sigma,\Omega=1,2,
\end{equation}
the torsion coefficients of the ON-section $\{N_1,N_2\}$ (or the coefficients of the normal connection of $X$).
\begin{proposition}
(The Weingarten equations)\\[0.1cm]
In $B$ there hold (use the summation convention for $j,$ $k,$ and $\Omega$)
\begin{equation}\label{5.3}
  N_{\Sigma,u^i}=-L_{\Sigma,ij}g^{jk}X_{u^k}+\sigma_{\Sigma,i}^\Omega N_\Omega\,,
  \quad i=1,2,\ \Sigma=1,2,
\end{equation}
where the $g^{ij}$ are defined via the Kronecker symbol $\delta_i^k$ as $g_{ij}g^{jk}=\delta_i^k.$
\end{proposition}
\noindent
The proof is essentially the same as for the Weingarten equations in $\mathbb R^3.$\\[1ex]
Various integrability conditions are associated with (\ref{5.3}). For example, if we evaluate the normal components of $N_{\Sigma,uv}-N_{\Sigma,vu}\equiv 0,$ we arrive at the
\begin{proposition}
(The Ricci equations)\\
In $B$ there hold the integrability conditions
\begin{equation}\label{5.4}
  \sigma_{\Sigma,1,v}^\Omega
  -\sigma_{\Sigma,2,u}^\Omega
  +\sigma_{\Sigma,1}^\Theta\sigma_{\Theta,2}^\Omega
  -\sigma_{\Sigma,2}^\Theta\sigma_{\Theta,1}^\Omega
  =(L_{\Sigma,1j}L_{\Omega,k2}-L_{\Sigma,2j}L_{\Omega,k1})g^{jk}
\end{equation}
for $\Sigma,\Omega=1,2.$
\end{proposition}
\begin{definition}
(The curvature tensor of the normal bundle)\\[0.1cm]
We define the curvature tensor of the normal bundle as
\begin{equation}\label{5.5}
  S_{\Sigma,ij}^\Omega
  :=\sigma_{\Sigma,i,u^j}^\Omega
    -\sigma_{\Sigma,j,u^i}^\Omega
    +\sigma_{\Sigma,i}^\Theta\sigma_{\Theta,j}^\Omega
    -\sigma_{\Sigma,j}^\Theta\sigma_{\Theta,i}^\Omega
\end{equation}
for $i,j=1,2$ and $\Sigma,\Omega=1,2.$
\end{definition}
\begin{remark}
The Ricci equations (\ref{5.4}) can be written in the form
\begin{equation}\label{5.6}
  S_{\Sigma,12}^\Omega=(L_{\Sigma,1j}L_{\Omega,k2}-L_{\Sigma,2j}L_{\Omega,k1})g^{jk}\,,\quad\Sigma,\Omega=1,2.
\end{equation}
\end{remark}
\subsubsection{Torsion-free ON-normal sections}
Let us construct a new ON-normal section $\{\widetilde N_1,\widetilde N_2\}$ via
\begin{equation}\label{5.7}
  \widetilde N_1=\cos\varphi N_1+\sin\varphi N_2\,,\quad
  \widetilde N_2=\sin\varphi N_1-\cos\varphi N_2
\end{equation}
with a rotation angle $\varphi.$ The associated new torsion coefficients $\widetilde\sigma_{\Sigma,i}^\Omega$ are
\begin{equation}\label{5.8}
  \widetilde\sigma_{1,1}^2=\sigma_{1,1}^2+\varphi_u\,,\quad
  \widetilde\sigma_{1,2}^2=\sigma_{1,2}^2+\varphi_v\,.
\end{equation}
\begin{problem}
Can we arrange $\{\widetilde N_1,\widetilde N_2\}$ such that $(\widetilde \sigma_{1,1}^2,\widetilde\sigma_{1,2}^2)\equiv(0,0)?$
\end{problem}
\noindent
Due to (\ref{5.8}) we have to solve the linear system
\begin{equation}\label{5.9}
  \varphi_u=-\sigma_{1,1}^2\,,\quad
  \varphi_v=-\sigma_{1,2}^2\,.
\end{equation}
At least locally, necessary and sufficient for the solvability of (\ref{5.9}) are the integrability conditions
\begin{equation}\label{5.10}
  0=\varphi_{uv}-\varphi_{vu}
   =-\sigma_{1,1,v}^2+\sigma_{1,2,u}^2\,.
\end{equation}
\begin{theorem}
(Flat normal bundles and torsion-free ON-normal sections)\\[0.1cm]
The vanishing of the curvature tensor (\ref{5.6}) is equivalent to the integrability conditions (\ref{5.10}).
\end{theorem}
\begin{proof}
Using $\sigma_{\Sigma,i}^\Sigma\equiv 0$ and $\sigma_{\Sigma,i}^\Omega=-\sigma_{\Omega,i}^\Sigma$ we calculate
\begin{itemize}
\item[$\bullet$]
$S_{1,ij}^1
 =\sigma_{1,i,u^j}^1-\sigma_{1,j,u^i}^1
  +\sigma_{1,i}^\Theta\sigma_{\Theta,j}^1-\sigma_{1,j}^\Theta\sigma_{\Theta,i}^1
 =\sigma_{1,i}^2\sigma_{2,j}^1-\sigma_{1,j}^2\sigma_{2,i}^1
 =0,$
\item[$\bullet$]
$S_{1,ij}^2
 =\sigma_{1,i,u^j}^2-\sigma_{1,j,u^i}^2
  +\sigma_{1,i}^\Theta\sigma_{\Theta,j}^2-\sigma_{1,j}^\Theta\sigma_{\Theta,i}^2
 =\sigma_{1,i,u^j}^2-\sigma_{1,j,u^i}^2\,,$
\item[$\bullet$]
$S_{2,ij}^1
 =\sigma_{2,i,u^j}^1-\sigma_{2,j,u^i}^1
  +\sigma_{2,i}^\Theta\sigma_{\Theta,j}^1-\sigma_{2,j}^\Theta\sigma_{\Theta,i}^1
 =\sigma_{2,i,u^j}^1-\sigma_{2,j,u^i}^1\,,$
\item[$\bullet$]
$S_{2,ij}^2
 =\sigma_{2,i,u^j}^2-\sigma_{2,j,u^i}^2
  +\sigma_{2,i}^\Theta\sigma_{\Theta,j}^2-\sigma_{2,j}^\Theta\sigma_{\Theta,i}^2
 =\sigma_{2,i}^1\sigma_{1,j}^2-\sigma_{2,j}^1\sigma_{1,i}^2
 =0.$
\end{itemize}
For the non-vanishing components of $S_{1,ij}^2$ and $S_{2,ij}^1$ we calculate
\begin{equation}\label{5.11}
\begin{array}{l}
  \displaystyle
  S_{1,11}^2=0,\quad S_{1,22}^2=0,\quad
  S_{1,12}^2=-S_{1,21}^2=\frac{\partial}{\partial v}\,\sigma_{1,1}^2-\frac{\partial}{\partial u}\,\sigma_{1,2}^2\,, \\[0.6cm]
  \displaystyle
  S_{2,11}^1=0,\quad S_{2,22}^1=0,\quad
  S_{2,12}^1=-S_{2,21}^1=\frac{\partial}{\partial v}\,\sigma_{2,1}^1-\frac{\partial}{\partial u}\,\sigma_{2,2}^1\,.
\end{array}
\end{equation}
These are exactly the conditions (\ref{5.10}).
\end{proof}
\subsection{Minimal surface graphs in $\mathbb R^4$}
The equivalence stated in the last theorem justifies the hope that curvature estimates and related Bernstein theorems for minimal graphs with a torsion-free ON-normal section $\{N_1,N_2\}$ are possible in the context of the theory of Erhard Heinz.\\[1ex]
We refer the reader to \cite{Wang_01} for results within the frame of mean curvature flow which combines curvature estimates and Bernstein type theorems for higher codimensional minimal graphs and certain stability criteria for the  second variation of the area functional.
\begin{theorem}
(Mu-Tao Wang, 2004)\\[0.1cm]
Let the complete minimal surface graph $\Sigma\subset\mathbb R^4$ be given with the properties:
\begin{itemize}
\item[(P1)]
$\mbox{\rm Area}\,(\Sigma\cap K_R)\le d_0R^2$ with some constant $d_0>0,$ and $K_R$ is the ball of radius $R$ in $\mathbb R^4$ centered at the origin;
\item[(P2)]
its normal bundle is flat, that is $S_{\Sigma,ij}^\Omega\equiv 0.$
\end{itemize}
Then $\Sigma$ is a plane.
\end{theorem}
\noindent
Because in our context up to now a detailed proof of this result fails we conclude our small note with the following remarks:
\begin{remarks}
\begin{itemize}
\item[1.]
For example, the minimal graph $(w,w^2)$ does not possess a torsion-free ON-normal section, that is its normal bundle is not flat.
\item[2.]
Immersions $X\colon B\to S^3$ have a flat normal bundle. But note that $\mathbb R^4$-minimal surfaces are not in $S^3.$
\item[3.]
Due to \cite{Wang_01}, Theorem 1.1, a minimal surface graph in $\mathbb R^4$ with flat normal bundle is stable.
\item[4.]
For stable immersions we can prove an estimate of the form (\ref{4.5}) for inner geodesic balls which could replace assumption (P1) eventually.
\item[5.]
We would like to draw your attention to \cite{SWX_01}, \cite{SWX_02}, where the authors study submanifolds with arbitrary codimension and prove Bernstein type theorems for minimal submanifolds with flat normal bundle.
\end{itemize}
\end{remarks}
\vspace*{0.4cm}
\noindent
Acknowledgement:\\[0.1cm]
The author wants to thank Matthias Bergner for many fruitful discussions.
\newpage\noindent
{\small

\vspace*{0.8cm}
\noindent
Steffen Fr\"ohlich\\
Technische Universit\"at Darmstadt\\
Fachbereich Mathematik, AG 4\\
Schlo\ss{}gartenstra\ss{}e 7\\
D-64289 Darmstadt\\
Germany\\[0.2cm]
e-mail: sfroehlich@mathematik.tu-darmstadt.de
}

\end{document}